\documentclass[11pt,a4paper]{article}

\usepackage[T1]{fontenc}
\usepackage[utf8]{inputenc}
\usepackage{amsmath,amssymb,amsthm,mathtools}
\usepackage[margin=1.1in]{geometry}
\usepackage{enumitem}
\usepackage[hidelinks]{hyperref}
\usepackage{tikz}
\usetikzlibrary{positioning,calc}

\theoremstyle{plain}
\newtheorem{theorem}{Theorem}[section]
\newtheorem{lemma}[theorem]{Lemma}
\newtheorem{proposition}[theorem]{Proposition}
\newtheorem{corollary}[theorem]{Corollary}

\theoremstyle{definition}
\newtheorem{definition}[theorem]{Definition}
\newtheorem{example}[theorem]{Example}
\newtheorem{assumption}[theorem]{Assumption}

\theoremstyle{remark}
\newtheorem{remark}[theorem]{Remark}

\newcommand{\R}{\mathbb{R}}
\newcommand{\N}{\mathbb{N}}
\newcommand{\Gt}{\widetilde{G}}
\newcommand{\Ea}{E_{\alpha}}
\newcommand{\Ga}{G_{\alpha}}
\newcommand{\Gta}{\widetilde{G}_{\alpha}}
\newcommand{\Na}{N(\alpha)}
\newcommand{\pa}{\varphi_{\alpha}}

\begin{document}

\title{Graded-Edge Rakotch Contractions in Generalized Metric Spaces
with Weighted Graphs}

\author{Ismat Beg\thanks{Department of Mathematics \& Statistical Sciences,
Lahore School of Economics, Lahore, Pakistan.
E-mail: \texttt{ibeg@lahoreschool.edu.pk}}
\and
Muhammad Aslam Javed\thanks{Department of Mathematics, Faculty of
Science \& Technology, University of Central Punjab, Lahore, Pakistan;
and Research Fellow, INTI International University, Persiaran Perdana BBN,
Putra Nilai, Malaysia.
E-mail: \texttt{m.a.javed@ucp.edu.pk}}}
\date{}

\maketitle

\begin{abstract}
We introduce the notion of graded-edge Rakotch contraction mappings on a generalized metric space
with a weighted graph. We prove that the successive iterate of every starting point converge uniformly to a common limit, which is the unique fixed point. We derive an explicit iteration-complexity bound and, under a relaxed triangle inequality, a geometric rate whose constant is the graph invariant. Finally, we illustrate the theory on an explicit real graded contraction that realizes the interior optimum, on a quadratic integral equation, and on finite weighted digraphs.
\end{abstract}



\section{Introduction}

Banach's contraction principle \cite{Banach1922} has led to a vast amount of literature, of which we mention only the parts directly relevant here. Rakotch
\cite{Rakotch1962} replaced the constant Lipschitz factor by a decreasing
comparison function $\varphi$ subject to $\varphi(t)<1$ for $t>0$, and thereby
obtained a strictly larger class of mappings retaining the conclusion of
Banach's theorem; the point of the relaxation is precisely that $\varphi(0)$ is
allowed to equal $1$, so that a Rakotch contraction need not be Lipschitz with a
constant smaller than one. A second line of development concerns spaces whose
distance function satisfies only a relaxed triangle inequality
$d(x,z)\le M\bigl(d(x,y)+d(y,z)\bigr)$ for a constant $M\ge1$. Such spaces, the
\emph{$b$-metric spaces} of Bakhtin \cite{Bakhtin1989} and Czerwik
\cite{Czerwik1993}, form a proper generalization of metric spaces, the case
$M=1$, and support their own extensive fixed point theory. A standard
instance is the snowflake $d(x,y)=|x-y|^{p}$ with $p>1$ on the real line, for
which $M=2^{p-1}$ and the ordinary triangle inequality fails.

Motivated by such examples, Reich and Zaslavski
\cite{RZ2022} introduced \emph{generalized metric spaces}, in which the
inequality $d(x,z)\le M d(x,y)+M d(y,z)$ is taken as an axiom, and proved a
Rakotch-type fixed point theorem. In the recent paper \cite{RZ2026} they
went further in two directions: First, they replaced the relaxed
triangle inequality by a purely qualitative $\varepsilon$--$\delta$ condition,
so that no numerical constant is available at all. Second, they endowed the
space with a graph $G$ and required the contractive estimate to hold only along
the edges of $G$, the two points of an arbitrary pair being connected instead by
a $G$-chain of length at most a fixed integer $n_{0}$. This last device goes
back to Jachymski \cite{Jachymski2008}, whose theorem unified the classical
contraction principle with its ordered-space variants due to Ran and Reurings
\cite{RanReurings2004} and Nieto and Rodr\'iguez-L\'opez \cite{NietoRL2005};
Beg et al. \cite{Beg2010} see also Bojor \cite{Bojor2012} and Beg and Butt
\cite{Beg2013}.

In \cite{RZ2026} the graph enters twice: once through the edge-preservation
requirement, and once through the connectivity hypothesis (P3). The contractive
estimate itself, however, is governed by a single scalar comparison function
$\varphi$, so that the graph contributes no quantitative information beyond the
combinatorial bound $n_{0}$. We propose to let the graph carry the contraction
data as well. Concretely, we attach to each edge a grade $w(e)\in(0,1]$, we let
$\Ea=\{e\in E(G):w(e)\ge\alpha\}$ denote the associated level cuts, and we ask
that the mapping contract with comparison function $\pa$ on $\Ea$, where
$\alpha\mapsto\pa$ is pointwise non-increasing.

Substituting $\varphi:=\varphi_{\alpha_{0}}$ with $\alpha_{0}=\inf_{e}w(e)$ restores the
ungraded hypothesis, so that the qualitative conclusions cannot change. What
makes the refinement substantive is that the two roles of the graph now move in
\emph{opposite} directions. Raising $\alpha$ improves the contraction factor but
deletes edges, so that the chain-length invariant
$\Na:=\sup_{x,y\in K}\rho_{\Gta}(x,y)$ increases. Every convergence estimate we
prove therefore depends on $\alpha$ through two competing terms and must be
optimized over the threshold. In Example~\ref{ex:interior} we exhibit an
explicit space, mapping and grading for which the optimal threshold is interior,
and for which the ungraded choice $\alpha_{0}$ is strictly worse than the graded
optimum by a factor that persists as the tolerance tends to zero. We believe this example and not the definition itself as justification of the framework.

The second observation is related to the origin of the grading. It does not need to be provided as additional information: it is already contained in any edge-contractive mapping. Defining $c(x,y)$ to be the supremum of the orbital ratios
$d(T^{n+1}x,T^{n+1}y)/d(T^{n}x,T^{n}y)$ and putting $w:=1-c$ produces a grade
which automatically satisfies the monotonicity axiom, and which normalizes the
comparison family to $\pa\equiv1-\alpha$ (Proposition~\ref{prop:canonical}).
The graded theory is therefore not an additional hypothesis imposed on $T$, but
a decomposition of the information $T$ already carries.

Our contributions are the following.
\begin{enumerate}[label=(\roman*),leftmargin=2.2em]
\item We introduce graded-edge Rakotch contractions
(Definitions~\ref{def:C1} and \ref{def:C2}) over the qualitative framework of
\cite{RZ2026}, retaining the possible asymmetry of $d$ and treating $G$ as a
directed graph throughout, and we show that the grading always exists
canonically (Proposition~\ref{prop:canonical}).
\item We make the chaining step of \cite{RZ2026} constructive by exhibiting an
explicit modulus $\eta_{m}$ obtained by iterating the modulus of the
$\varepsilon$--$\delta$ triangle condition (Lemma~\ref{lem:chain}).
\item We prove a global uniform convergence theorem
(Theorem~\ref{thm:A}) which contains \cite[Theorem 2]{RZ2026} as the case
$w\equiv1$, and we close a gap in \cite[Theorem 2]{RZ2026}, whose statement
asserts the convergence of the Picard iterates but not the fixed point property
of the limit (Theorem~\ref{thm:B}).
\item We derive an explicit iteration-complexity bound
(Theorem~\ref{thm:C}) and, under a relaxed triangle inequality, a geometric rate
whose constant is the graph invariant $\Na M^{\Na-1}$
(Theorem~\ref{thm:D}); we optimize over $\alpha$ (Corollary~\ref{cor:alphaopt})
and settle the resulting trade-off explicitly (Example~\ref{ex:interior}).
\item We illustrate the theory on an explicit real graded contraction that
realizes the interior optimum, on a quadratic integral equation, and on finite
weighted digraphs (Section~\ref{sec:examples}).
\end{enumerate}

\section{Preliminaries}

\subsection{The space}

Throughout, $X$ is a nonempty set and $d\colon X\times X\to[0,\infty)$ satisfies
\begin{equation}\label{eq:d0}
d(x,y)=0 \iff x=y \qquad (x,y\in X).
\end{equation}
Following \cite{RZ2026} we do \emph{not} assume that $d$ is symmetric, nor that
it satisfies a triangle inequality; instead we impose the following two
conditions.

\begin{assumption}[A1]\label{ass:A1}
For each $\varepsilon>0$ there exists $\delta>0$ such that
$d(x,y)\le\delta$ implies $d(y,x)\le\varepsilon$.
\end{assumption}

\begin{assumption}[A2]\label{ass:A2}
For each $\varepsilon>0$ there exists $\delta>0$ such that
$d(x,y)\le\delta$ and $d(y,z)\le\delta$ imply $d(x,z)\le\varepsilon$.
\end{assumption}

For $\varepsilon>0$ we fix once and for all a choice of 
$\delta_{1}(\varepsilon)$ for \textup{(A1)} and $\delta_{2}(\varepsilon)$ for
\textup{(A2)}. Without loss of generality, we may  assume that
\begin{equation}\label{eq:moduli}
\delta_{i}\ \text{is non-decreasing and}\quad
\delta_{i}(\varepsilon)\le\varepsilon \qquad (i=1,2,\ \varepsilon>0).
\end{equation}
We refer to $\delta_{1}$ and $\delta_{2}$ as the \emph{symmetry modulus} and the
\emph{triangle modulus} of $(X,d)$.

Setting
\begin{equation}\label{eq:unif}
U(\varepsilon):=\{(x,y)\in X\times X: d(x,y)\le\varepsilon,\ d(y,x)\le\varepsilon\},
\qquad \varepsilon>0,
\end{equation}
one checks as in \cite{RZ2026} that $\{U(\varepsilon)\}_{\varepsilon>0}$ is a base
of a uniformity on $X$, metrizable by a metric $d_{1}$. We assume throughout
that the metric space $(X,d_{1})$ is complete. By construction, for every
$\rho>0$ there is $\varepsilon>0$ with
$U(\varepsilon)\subseteq\{(x,y):d_{1}(x,y)\le\rho\}$, and conversely for every
$\varepsilon>0$ there is $\rho>0$ with
$\{(x,y):d_{1}(x,y)\le\rho\}\subseteq U(\varepsilon)$; in particular
\begin{equation}\label{eq:d1d}
d_{1}(x_{n},x)\to0 \ \Longrightarrow\ d(x_{n},x)\to0 \ \text{ and }\
d(x,x_{n})\to0 .
\end{equation}

A set $E\subseteq X$ is \emph{bounded} if
$\sup\{d(x,y):x,y\in E\}<\infty$.

\begin{remark}\label{rem:regimeB}
In several places below we shall strengthen \textup{(A2)} to the quantitative
condition
\begin{equation}\label{eq:relaxed}
d(x,z)\le M\bigl(d(x,y)+d(y,z)\bigr)\qquad (x,y,z\in X),
\end{equation}
which is the notion of a generalized metric space of \cite{RZ2022}; when $d$ is
symmetric this is exactly a $b$-metric space with constant $M$ in the sense of
\cite{Bakhtin1989,Czerwik1993}. Taking $y=x$ in \eqref{eq:relaxed} forces
$M\ge1$. Condition \eqref{eq:relaxed} implies
\textup{(A2)} with $\delta_{2}(\varepsilon)=\varepsilon/(2M)$. We refer to the
setting of \textup{(A1)}--\textup{(A2)} as the \emph{qualitative regime} and to
the setting of \eqref{eq:relaxed} as the \emph{quantitative regime}. When $M=1$
and $d$ is symmetric, $(X,d)$ is an ordinary metric space.
\end{remark}

\subsection{Comparison functions}

Let $\Phi$ denote the set of all decreasing functions
$\varphi\colon[0,\infty)\to[0,1]$ such that $\varphi(t)<1$ for every $t>0$.

\begin{remark}\label{rem:codomain}
The requirement $\varphi(t)<1$ is imposed only for $t>0$, so that
$\varphi(0)=1$ is permitted. This is essential: a decreasing
$\varphi\colon[0,\infty)\to[0,1)$ would satisfy $\varphi(t)\le\varphi(0)<1$ for
every $t$, so that the Rakotch condition would collapse to the Banach condition
with constant $\varphi(0)$ and no generality would be gained. The definition in
\cite[\S1]{RZ2026} has codomain $[0,1]$, as here; the one in \cite[\S2]{RZ2026}
has codomain $[0,1)$, which we take to be a slip of the pen.
\end{remark}

\subsection{Graded graphs}

Let $G$ be a directed graph with vertex set $V(G)=X$ and edge set
$E(G)\subseteq X\times X$ containing all loops. We write $\Gt$ for the
symmetrization of $G$, that is, the undirected graph in which $x$ and $y$ are
adjacent whenever $(x,y)\in E(G)$ or $(y,x)\in E(G)$.

\begin{definition}\label{def:graded}
A \emph{grade} on $G$ is a map $w\colon E(G)\to(0,1]$. For $\alpha\in(0,1]$ the
\emph{level cut} of $G$ at $\alpha$ is the directed graph $\Ga$ with
$V(\Ga)=X$ and
\[
\Ea:=\{e\in E(G): w(e)\ge\alpha\}.
\]
The pair $(G,w)$ is called a \emph{graded graph} on $X$.
\end{definition}

Clearly $E_{1}\subseteq\Ea\subseteq E_{\beta}$ whenever $\beta\le\alpha$, so the
level cuts form a decreasing family of graphs.

Given a nonempty $K\subseteq X$, a \emph{$\Gta$-chain} in $K$ from $x$ to $y$ of
length $q$ is a finite sequence $\xi_{0},\dots,\xi_{q}\in K$ with $\xi_{0}=x$,
$\xi_{q}=y$ and, for each $i\in\{0,\dots,q-1\}$,
\[
(\xi_{i},\xi_{i+1})\in\Ea \quad\text{or}\quad (\xi_{i+1},\xi_{i})\in\Ea .
\]
We denote by $\rho_{\Gta}(x,y)$ the least length of such a chain (with
$\rho_{\Gta}(x,y)=\infty$ if none exists) and we set
\begin{equation}\label{eq:Nalpha}
\Na:=\sup_{x,y\in K}\rho_{\Gta}(x,y)\in\N\cup\{\infty\},
\end{equation}
the \emph{$\alpha$-chain diameter} of $K$. The function $\alpha\mapsto\Na$ is
non-decreasing, and
\[
A:=\{\alpha\in(0,1]: \Na<\infty\}
\]
is a sub-interval of $(0,1]$ containing $0$ in its closure whenever it is
nonempty. The connectivity hypothesis (P3) of \cite[page 13]{RZ2026} is precisely the statement
$\Na<\infty$ for the ungraded graph, with $n_{0}=\Na$.

\section{Graded-edge Rakotch contractions}\label{sec:defs}

Throughout this section, $(G,w)$ is a graded graph on $X$, $K\subseteq X$ is
nonempty and closed in $(X,d_{1})$, and $T\colon K\to K$.

\begin{definition}[Grade-monotone edge preservation, (C1)]\label{def:C1}
The mapping $T$ is said to be \emph{grade-monotone edge preserving} if for every
$(x,y)\in E(G)$ with $x,y\in K$,
\[
(T(x),T(y))\in E(G)\qquad\text{and}\qquad w\bigl(T(x),T(y)\bigr)\ \ge\ w(x,y).
\]
\end{definition}

Condition (C1) strengthens the edge-preservation requirement of
\cite[(5)]{RZ2026} by demanding that $T$ never weakens a connection. It is what
guarantees that an orbit which starts inside a level cut remains there.

\begin{definition}[Graded Rakotch condition, (C2)]\label{def:C2}
Let $\{\pa\}_{\alpha\in(0,1]}\subseteq\Phi$ be a family such that
\begin{equation}\label{eq:monotonefamily}
\alpha\le\beta \ \Longrightarrow\ \varphi_{\beta}(t)\le\pa(t)
\qquad (t\ge0),
\end{equation}
that is, $\alpha\mapsto\pa$ is pointwise non-increasing. The mapping $T$ is a
\emph{graded-edge Rakotch contraction} with respect to $(G,w)$ and
$\{\pa\}$ if it satisfies \textup{(C1)} and
\begin{equation}\label{eq:C2}
d\bigl(T(x),T(y)\bigr)\ \le\ \varphi_{w(x,y)}\bigl(d(x,y)\bigr)\, d(x,y)
\qquad\text{for every }(x,y)\in E(G)\text{ with }x,y\in K.
\end{equation}
\end{definition}

We set $T^{0}(x)=x$ for $x\in K$.

\begin{remark}\label{rem:reduction}
If $w\equiv1$ and $\varphi_{1}=\varphi$, then \textup{(C1)} reduces to
\cite[(5)]{RZ2026} and \eqref{eq:C2} reduces to \cite[(6)]{RZ2026}. If moreover
$E(G)=K\times K$, one recovers the hypotheses of \cite[Theorem 1]{RZ2026}.
Finally, if each $\pa$ is constant, one recovers the graph-theoretic Banach
contractions of \cite{Jachymski2008}.
\end{remark}

The next result shows that Definition~\ref{def:C2} does not in fact impose any
new hypothesis on $T$ beyond edge preservation and a mild uniformity: the
grading may always be manufactured from the mapping itself, and when it is, the
comparison family is normalized.

\begin{proposition}[Canonical grading]\label{prop:canonical}

Let $T$ satisfy $(x,y)\in E(G)\Rightarrow (T(x),T(y))\in E(G)$ for $x,y\in K$.
For $(x,y)\in E(G)$ with $x,y\in K$ put
\begin{equation}\label{eq:cdef}
c(x,y):=\sup\left\{
\frac{d\bigl(T^{n+1}(x),T^{n+1}(y)\bigr)}{d\bigl(T^{n}(x),T^{n}(y)\bigr)}
\ :\ n\ge0,\ d\bigl(T^{n}(x),T^{n}(y)\bigr)>0\right\},
\end{equation}
with the convention $\sup\emptyset:=0$, and assume $c(x,y)<1$ for every such
pair. Then $w:=1-c$ is a grade on $E(G)$, and $T$ is a graded-edge Rakotch
contraction with respect to $(G,w)$ and the family
\[
\pa\equiv 1-\alpha \qquad (\alpha\in(0,1]).
\]
Moreover $1-\alpha$ is optimal, in the sense that no constant $c'<1-\alpha$
satisfies $d(T(x),T(y))\le c'\,d(x,y)$ for every $(x,y)\in\Ea$, whenever the
value $\alpha$ is attained by $w$ and the supremum in \eqref{eq:cdef} is
attained at $n=0$ for some such edge.
\end{proposition}

\begin{proof}
By hypothesis $c(x,y)\in[0,1)$, so $w(x,y)=1-c(x,y)\in(0,1]$ and $w$ is a grade.

\emph{(C1).} Edge preservation holds by assumption. Since the mapping $T$ preserves edges,
the pair $(T(x),T(y))$ lies in $E(G)$ and, directly from \eqref{eq:cdef},
\[
c\bigl(T(x),T(y)\bigr)
=\sup\left\{
\frac{d\bigl(T^{n+1}(x),T^{n+1}(y)\bigr)}{d\bigl(T^{n}(x),T^{n}(y)\bigr)}
: n\ge1,\ d\bigl(T^{n}(x),T^{n}(y)\bigr)>0\right\}
\ \le\ c(x,y),
\]
the supremum on the right being taken over a subset of the index set defining
$c(x,y)$. Hence $w(T(x),T(y))\ge w(x,y)$.

\emph{(C2).} Each $\pa\equiv1-\alpha$ is a constant function with value in
$[0,1)$, hence lies in $\Phi$, and $\alpha\mapsto1-\alpha$ is decreasing, so
\eqref{eq:monotonefamily} holds. Let $(x,y)\in E(G)$ with $x,y\in K$. If
$d(x,y)=0$ then $x=y$ and \eqref{eq:C2} is trivial. Otherwise the case $n=0$ of
\eqref{eq:cdef} gives
\[
d\bigl(T(x),T(y)\bigr)\le c(x,y)\,d(x,y)=\bigl(1-w(x,y)\bigr)d(x,y)
=\varphi_{w(x,y)}\bigl(d(x,y)\bigr)\,d(x,y).
\]

\emph{Optimality.} Suppose $(x,y)\in\Ea$ with $w(x,y)=\alpha$ and that the
supremum in \eqref{eq:cdef} is attained at $n=0$. Then
$d(T(x),T(y))=c(x,y)d(x,y)=(1-\alpha)d(x,y)$, so no smaller constant is
admissible.
\end{proof}

\begin{remark}\label{rem:notcosmetic}
It is tempting to dismiss Definition~\ref{def:C2} on the ground that
$\varphi:=\varphi_{\alpha_{0}}$ with
$\alpha_{0}:=\inf\{w(e):e\in E(G)\}$ satisfies \eqref{eq:C2} on all of $E(G)$,
so that every graded-edge Rakotch contraction is an ordinary one. This is
correct, and it is exactly why the qualitative statements of
Section~\ref{sec:main} are not new in substance. The content of the grading lies
in Section~\ref{sec:quant}: the quantitative estimates involve $\pa$ and $\Na$
simultaneously, and passing to $\alpha_{0}$ minimizes the first factor while
maximizing the second. Example~\ref{ex:interior} exhibits data for which the
choice $\alpha_{0}$ is strictly worse than the optimal graded choice.
\end{remark}

\section{Auxiliary results}

We first record the effect of \textup{(C1)} along an orbit.

\begin{lemma}[Invariance and monotonicity]\label{lem:mono}
Let $T$ be a graded-edge Rakotch contraction, let $\alpha\in(0,1]$ and let
$x,y\in K$ satisfy $(x,y)\in\Ea$. Then
\[
\bigl(T^{n}(x),T^{n}(y)\bigr)\in\Ea \qquad (n\ge0),
\]
and the sequence $n\mapsto d\bigl(T^{n}(x),T^{n}(y)\bigr)$ is non-increasing.
If moreover $\pa\equiv c_{\alpha}$ is constant, then
\begin{equation}\label{eq:geom}
d\bigl(T^{n}(x),T^{n}(y)\bigr)\ \le\ c_{\alpha}^{\,n}\,d(x,y)\qquad(n\ge0).
\end{equation}
\end{lemma}

\begin{proof}
The first assertion holds for $n=0$ by hypothesis. If
$(T^{n}(x),T^{n}(y))\in\Ea$, then \textup{(C1)} gives
$(T^{n+1}(x),T^{n+1}(y))\in E(G)$ and
\[
w\bigl(T^{n+1}(x),T^{n+1}(y)\bigr)\ \ge\ w\bigl(T^{n}(x),T^{n}(y)\bigr)\ \ge\ \alpha,
\]
so that $(T^{n+1}(x),T^{n+1}(y))\in\Ea$; the claim follows by induction.
For the second assertion, put $a_{n}:=d(T^{n}(x),T^{n}(y))$. Applying
\eqref{eq:C2} to the edge $(T^{n}(x),T^{n}(y))\in E(G)$ and using
$\varphi_{w}\le1$ we obtain
\[
a_{n+1}\ \le\ \varphi_{w(T^{n}(x),T^{n}(y))}(a_{n})\,a_{n}\ \le\ a_{n}.
\]
If $\pa\equiv c_{\alpha}$, then $w(T^{n}(x),T^{n}(y))\ge\alpha$ and
\eqref{eq:monotonefamily} give $a_{n+1}\le c_{\alpha}a_{n}$, whence
\eqref{eq:geom} by induction.
\end{proof}

The next lemma is the graded form of \cite[Lemma 1]{RZ2026}. It converts the
multiplicative contraction \eqref{eq:C2} into an additive decrease that is
uniform over all edges of a given level cut.

\begin{lemma}[Graded uniform decrease]\label{lem:decrease}
Let $T$ be a graded-edge Rakotch contraction, let $\alpha\in(0,1]$, let
$\varepsilon>0$, and let $x,y\in K$ satisfy
\[
(x,y)\in\Ea \qquad\text{and}\qquad d(x,y)>\varepsilon .
\]
Then
\[
d\bigl(T(x),T(y)\bigr)\ \le\ d(x,y)-\varepsilon\bigl(1-\pa(\varepsilon)\bigr).
\]
\end{lemma}

\begin{proof}
Write $w:=w(x,y)$ and $t:=d(x,y)$. Since $(x,y)\in\Ea$ we have $w\ge\alpha$, so
that $\varphi_{w}(t)\le\pa(t)$ by \eqref{eq:monotonefamily}. Since $\pa$ is
decreasing and $t>\varepsilon$, we also have $\pa(t)\le\pa(\varepsilon)$. Hence,
by \eqref{eq:C2},
\[
t-d\bigl(T(x),T(y)\bigr)\ \ge\ t\bigl(1-\varphi_{w}(t)\bigr)
\ \ge\ t\bigl(1-\pa(t)\bigr)
\ \ge\ t\bigl(1-\pa(\varepsilon)\bigr)
\ \ge\ \varepsilon\bigl(1-\pa(\varepsilon)\bigr),
\]
where the last inequality uses $t>\varepsilon$ together with
$1-\pa(\varepsilon)\ge0$. Rearranging gives the assertion.
\end{proof}

In \cite{RZ2026} the existence of a threshold $\varepsilon_{0}$ with the property
that short steps along a bounded chain produce a short displacement between its
endpoints is deduced from \textup{(A2)} but is not exhibited. The following
lemma supplies it explicitly. Recall the moduli $\delta_{1},\delta_{2}$ fixed in
\eqref{eq:moduli} and let $\delta_{2}^{\circ k}$ denote the $k$-fold composition
of $\delta_{2}$ with itself, with $\delta_{2}^{\circ0}:=\mathrm{id}$.

\begin{lemma}[Chaining modulus]\label{lem:chain}
Define
\begin{equation}\label{eq:etadef}
\eta_{1}(\varepsilon):=\min\{\varepsilon,\delta_{1}(\varepsilon)\},
\qquad
\eta_{m}(\varepsilon):=\eta_{1}\bigl(\delta_{2}^{\circ(m-1)}(\varepsilon)\bigr)
\quad (m\ge2,\ \varepsilon>0).
\end{equation}
Let $m\ge1$, let $\varepsilon>0$ and let $\xi_{0},\dots,\xi_{m}\in X$ be such
that for each $i\in\{0,\dots,m-1\}$ at least one of the inequalities
\[
d(\xi_{i},\xi_{i+1})\le\eta_{m}(\varepsilon),
\qquad
d(\xi_{i+1},\xi_{i})\le\eta_{m}(\varepsilon)
\]
holds. Then $d(\xi_{0},\xi_{m})\le\varepsilon$.
\end{lemma}

\begin{proof}
We argue by induction on $m$. Observe first that, by \eqref{eq:moduli}, each
$\eta_{m}$ is non-decreasing and $\eta_{m+1}\le\eta_{m}\le\eta_{1}\le\mathrm{id}$.

\emph{Case $m=1$.} If $d(\xi_{0},\xi_{1})\le\eta_{1}(\varepsilon)\le\varepsilon$
there is nothing to prove. If instead
$d(\xi_{1},\xi_{0})\le\eta_{1}(\varepsilon)\le\delta_{1}(\varepsilon)$, then
Assumption~\ref{ass:A1} yields $d(\xi_{0},\xi_{1})\le\varepsilon$.

\emph{Inductive step.} Let $m\ge1$ and assume the assertion holds for chains of
length $m$. Let $\xi_{0},\dots,\xi_{m+1}$ be a chain all of whose steps are
bounded by $\eta_{m+1}(\varepsilon)$, and put
$\varepsilon':=\delta_{2}(\varepsilon)$. By \eqref{eq:etadef},
\[
\eta_{m+1}(\varepsilon)
=\eta_{1}\bigl(\delta_{2}^{\circ m}(\varepsilon)\bigr)
=\eta_{1}\bigl(\delta_{2}^{\circ(m-1)}(\varepsilon')\bigr)
=\eta_{m}(\varepsilon').
\]
Hence the first $m$ steps of the chain are bounded by $\eta_{m}(\varepsilon')$,
and the induction hypothesis, applied with $\varepsilon'$ in place of
$\varepsilon$, gives
\[
d(\xi_{0},\xi_{m})\ \le\ \varepsilon'=\delta_{2}(\varepsilon).
\]
The last step is bounded by
$\eta_{m+1}(\varepsilon)=\eta_{m}(\varepsilon')\le\eta_{1}(\varepsilon')$, so the
case $m=1$ applied to the pair $(\xi_{m},\xi_{m+1})$ with $\varepsilon'$ in place
of $\varepsilon$ gives
\[
d(\xi_{m},\xi_{m+1})\ \le\ \varepsilon'=\delta_{2}(\varepsilon).
\]
Assumption~\ref{ass:A2} now yields $d(\xi_{0},\xi_{m+1})\le\varepsilon$, which
completes the induction.
\end{proof}

\begin{remark}\label{rem:etaquant}
In the quantitative regime \eqref{eq:relaxed} with $d$ symmetric one may take
$\delta_{1}(\varepsilon)=\varepsilon$ and
$\delta_{2}(\varepsilon)=\varepsilon/(2M)$, whence
\[
\eta_{m}(\varepsilon)=\frac{\varepsilon}{(2M)^{m-1}}\qquad(m\ge1).
\]
Thus the chaining modulus degrades geometrically in the chain length, at a rate
governed by the constant $M$. This is the precise mechanism by which the graph
invariant $\Na$ enters the estimates of Section~\ref{sec:quant}, and the reason
why no analogue of those estimates is available in the purely qualitative
regime.
\end{remark}

\section{Fixed point}\label{sec:main}

Let $K\subseteq X$ be a nonempty and closed in $(X,d_{1})$,
\begin{equation}\label{eq:D0}
D_{0}:=\sup\{d(x,y):x,y\in K\}<\infty,
\end{equation}
$T\colon K\to K$ is a graded-edge Rakotch contraction with respect to a graded
graph $(G,w)$ and a family $\{\pa\}$ as in Definition~\ref{def:C2}, and
$A\ne\emptyset$, that is, $\Na<\infty$ for at least one $\alpha\in(0,1]$.

For $\varepsilon>0$ and $\alpha\in A$ we introduce the two quantities
\begin{equation}\label{eq:eps0k0}
\varepsilon_{0}(\varepsilon,\alpha):=\eta_{\Na}(\varepsilon),
\qquad
k_{0}(\varepsilon,\alpha):=
\left\lceil
\frac{D_{0}}
{\varepsilon_{0}(\varepsilon,\alpha)\,
\bigl(1-\pa\bigl(\varepsilon_{0}(\varepsilon,\alpha)\bigr)\bigr)}
\right\rceil .
\end{equation}

\begin{assumption}[A3]\label{ass:A3}
$\varepsilon_{0}(\varepsilon,\alpha)$ and $k_{0}(\varepsilon,\alpha)$ are finite and positive.
\end{assumption}

\begin{theorem}[Global uniform convergence]\label{thm:A}
Let Assumption~\ref{ass:A3} hold and let $\alpha\in A$. Then:
\begin{enumerate}[label=\textup{(\roman*)},leftmargin=2.4em]
\item for every $\varepsilon>0$ and every integer $n\ge k_{0}(\varepsilon,\alpha)$,
\begin{equation}\label{eq:unifsmall}
\sup_{x,y\in K} d\bigl(T^{n}(x),T^{n}(y)\bigr)\ \le\ \varepsilon;
\end{equation}
\item there exists a point $x^{*}\in K$ such that
$T^{n}(x)\to x^{*}$ in $(X,d_{1})$ as $n\to\infty$, for every $x\in K$, and the
convergence is uniform on $K$.
\end{enumerate}
\end{theorem}

\begin{proof}
Write $\varepsilon_{0}:=\varepsilon_{0}(\varepsilon,\alpha)$ and
$k_{0}:=k_{0}(\varepsilon,\alpha)$.

\emph{Step 1 (a single edge enters the $\varepsilon_{0}$-regime in at most
$k_{0}$ steps).} Let $(\xi,\eta)\in\Ea$ with $\xi,\eta\in K$. We claim that
there is an integer $n\in\{0,\dots,k_{0}\}$ with
$d(T^{n}(\xi),T^{n}(\eta))\le\varepsilon_{0}$. Suppose not. Then
$d(T^{n}(\xi),T^{n}(\eta))>\varepsilon_{0}$ for every
$n\in\{0,\dots,k_{0}\}$, and since $(T^{n}(\xi),T^{n}(\eta))\in\Ea$ by
Lemma~\ref{lem:mono}, Lemma~\ref{lem:decrease} yields
\[
d\bigl(T^{n+1}(\xi),T^{n+1}(\eta)\bigr)
\ \le\ d\bigl(T^{n}(\xi),T^{n}(\eta)\bigr)
-\varepsilon_{0}\bigl(1-\pa(\varepsilon_{0})\bigr),
\qquad n\in\{0,\dots,k_{0}-1\}.
\]
Summing this telescoping estimate and using \eqref{eq:D0} we obtain
\[
\begin{split}
D_{0}\ \ge\ d(\xi,\eta)-d\bigl(T^{k_{0}}(\xi),T^{k_{0}}(\eta)\bigr)
&=\sum_{i=0}^{k_{0}-1}
\Bigl[d\bigl(T^{i}(\xi),T^{i}(\eta)\bigr)-d\bigl(T^{i+1}(\xi),T^{i+1}(\eta)\bigr)\Bigr]\\
&\ \ge\ k_{0}\,\varepsilon_{0}\bigl(1-\pa(\varepsilon_{0})\bigr),
\end{split}
\]
that is, $k_{0}\le D_{0}\varepsilon_{0}^{-1}(1-\pa(\varepsilon_{0}))^{-1}$. This
contradicts \eqref{eq:eps0k0} and proves the claim.

\emph{Step 2 (the regime is absorbing).} By Lemma~\ref{lem:mono} the sequence
$n\mapsto d(T^{n}(\xi),T^{n}(\eta))$ is non-increasing. Combined with Step 1
this shows that
\begin{equation}\label{eq:edgesmall}
d\bigl(T^{n}(\xi),T^{n}(\eta)\bigr)\le\varepsilon_{0}
\qquad\text{for every }n\ge k_{0}\text{ and every }(\xi,\eta)\in\Ea
\text{ with }\xi,\eta\in K.
\end{equation}

\emph{Step 3 (chaining).} Let $x,y\in K$. Since $\Na<\infty$ there are an
integer $q\le\Na$ and points $\xi_{0},\dots,\xi_{q}\in K$ forming a
$\Gta$-chain from $x$ to $y$. Fix $n\ge k_{0}$. By \textup{(C1)} the points
$T^{n}(\xi_{0}),\dots,T^{n}(\xi_{q})$ again form a $\Gta$-chain, and by
\eqref{eq:edgesmall} each of its steps satisfies
\[
d\bigl(T^{n}(\xi_{i}),T^{n}(\xi_{i+1})\bigr)\le\varepsilon_{0}
\quad\text{or}\quad
d\bigl(T^{n}(\xi_{i+1}),T^{n}(\xi_{i})\bigr)\le\varepsilon_{0},
\qquad i\in\{0,\dots,q-1\}.
\]
Since $\varepsilon_{0}=\eta_{\Na}(\varepsilon)\le\eta_{q}(\varepsilon)$,
Lemma~\ref{lem:chain} applies and gives
$d(T^{n}(x),T^{n}(y))\le\varepsilon$. As $x,y\in K$ were arbitrary, this proves
\eqref{eq:unifsmall} and hence (i).

\emph{Step 4 (the orbits are Cauchy).} Let $\rho>0$ and choose $\varepsilon>0$
so small that $U(\varepsilon)\subseteq\{(u,v):d_{1}(u,v)\le\rho\}$; replacing
$\varepsilon$ by $\min\{\varepsilon,\delta_{1}(\varepsilon)\}$ we may assume in
addition that $d(u,v)\le\varepsilon$ implies $d(v,u)\le\varepsilon$, so that
\eqref{eq:unifsmall} yields
\begin{equation}\label{eq:inU}
\bigl(T^{n}(x),T^{n}(y)\bigr)\in U(\varepsilon)
\qquad (n\ge k_{0}(\varepsilon,\alpha),\ x,y\in K).
\end{equation}
Given $x\in K$ and $m\ge0$, apply \eqref{eq:inU} with $y:=T^{m}(x)\in K$ to get
$d_{1}(T^{n}(x),T^{n+m}(x))\le\rho$ for all $n\ge k_{0}(\varepsilon,\alpha)$.
Hence $\{T^{n}(x)\}_{n\ge1}$ is a Cauchy sequence in $(X,d_{1})$. Since
$(X,d_{1})$ is complete and $K$ is closed, the limit
$x^{*}_{x}:=\lim_{n\to\infty}T^{n}(x)$ exists and belongs to $K$.

\emph{Step 5 (the limit does not depend on the starting point).} Let
$x,y\in K$ and $\rho>0$. By \eqref{eq:inU}, $d_{1}(T^{n}(x),T^{n}(y))\le\rho$
for all sufficiently large $n$. Letting $n\to\infty$ and using the continuity of
$d_{1}$ we obtain $d_{1}(x^{*}_{x},x^{*}_{y})\le\rho$, and since $\rho>0$ was
arbitrary, $x^{*}_{x}=x^{*}_{y}=:x^{*}$. Finally, the integer
$k_{0}(\varepsilon,\alpha)$ of \eqref{eq:eps0k0} does not depend on $x$, so the
convergence is uniform on $K$. This proves (ii).
\end{proof}

\begin{corollary}\label{cor:RZ}
Taking $w\equiv1$, $\varphi_{1}=\varphi$ and $\alpha=1$ in
Theorem~\ref{thm:A} recovers \cite[Theorem 2]{RZ2026}, with $n_{0}=N(1)$.
\end{corollary}

Theorem~\ref{thm:A} asserts that all Picard orbits converge to a common limit,
but it does not by itself identify that limit as a fixed point of $T$; the same
is true of \cite[Theorem 2]{RZ2026}. The obstruction is that $T$ is assumed to
be contractive only along edges, so that no continuity of $T$ at $x^{*}$ is
available a priori. We record two hypotheses, either of which removes it. The
second is the graph-theoretic one and is modelled on property (A) of
\cite{Jachymski2008}.

\begin{definition}\label{def:propA}
The graded graph $(G,w)$ has \emph{property} $(\mathrm{A}_{\alpha})$ on $K$ if,
whenever $\{x_{n}\}_{n\ge1}\subseteq K$ satisfies $(x_{n},x_{n+1})\in\Ea$ for all
$n$ and $x_{n}\to x$ in $(X,d_{1})$ with $x\in K$, there exists a subsequence
$\{x_{n_{k}}\}$ with $(x_{n_{k}},x)\in\Ea$ for all $k$.
\end{definition}

\begin{theorem}[Fixed point and uniqueness]\label{thm:B}
Let Assumption~\ref{ass:A3} hold, let $\alpha\in A$, and let $x^{*}$ be as
in Theorem~\ref{thm:A}. Assume that at least one of the following conditions
holds:
\begin{enumerate}[label=\textup{(\alph*)},leftmargin=2.4em]
\item $T$ is continuous at $x^{*}$ as a mapping of $(K,d_{1})$ into itself;
\item $(G,w)$ has property $(\mathrm{A}_{\alpha})$ on $K$ and there exists
$x_{0}\in K$ with $(x_{0},T(x_{0}))\in\Ea$.
\end{enumerate}
Then $T(x^{*})=x^{*}$, and $x^{*}$ is the unique fixed point of $T$ in $K$.
\end{theorem}

\begin{proof}
Throughout, $x^{*}$ denotes the point produced by Theorem~\ref{thm:A}(ii); we
recall from that theorem that $x^{*}\in K$ and that $T^{n}(x)\to x^{*}$ in
$(X,d_{1})$ for \emph{every} $x\in K$. The proof falls into three parts: the
fixed point property under (a), the fixed point property under (b), and
uniqueness, the last of which uses neither (a) nor (b).

\medskip
\noindent\textbf{Part I: $T(x^{*})=x^{*}$ under hypothesis (a).}

\emph{Step 1.} Since $x^{*}\in K$, Theorem~\ref{thm:A}(ii) may be applied with
the starting point $x:=x^{*}$, and gives
\begin{equation}\label{eq:B1}
T^{n}(x^{*})\ \longrightarrow\ x^{*}\quad\text{in }(X,d_{1}).
\end{equation}

\emph{Step 2.} The points $T^{n}(x^{*})$ all lie in $K$, because $T$ maps $K$
into itself. Hypothesis (a) states that $T$ is continuous at $x^{*}$ as a
mapping of $(K,d_{1})$ into itself; applying this to the convergent sequence
\eqref{eq:B1} gives
\begin{equation}\label{eq:B2}
T\bigl(T^{n}(x^{*})\bigr)\ \longrightarrow\ T(x^{*})\quad\text{in }(X,d_{1}).
\end{equation}

\emph{Step 3.} On the other hand $T(T^{n}(x^{*}))=T^{n+1}(x^{*})$, and the
sequence $\{T^{n+1}(x^{*})\}_{n\ge0}$ is a shift of the sequence in
\eqref{eq:B1}; a shift of a convergent sequence converges to the same limit, so
\begin{equation}\label{eq:B3}
T\bigl(T^{n}(x^{*})\bigr)\ \longrightarrow\ x^{*}\quad\text{in }(X,d_{1}).
\end{equation}

\emph{Step 4.} $(X,d_{1})$ is a metric space, so a sequence has at most one
limit. Comparing \eqref{eq:B2} and \eqref{eq:B3} gives $T(x^{*})=x^{*}$.

\medskip
\noindent\textbf{Part II: $T(x^{*})=x^{*}$ under hypothesis (b).}

Let $x_{0}\in K$ be the point supplied by (b), so that
$(x_{0},T(x_{0}))\in\Ea$, and put $x_{n}:=T^{n}(x_{0})$ for $n\ge0$.

\emph{Step 1 (the orbit is a path in $\Ea$).} Apply Lemma~\ref{lem:mono} to the
pair $(x,y):=(x_{0},T(x_{0}))\in\Ea$. Its first assertion gives
$\bigl(T^{n}(x_{0}),T^{n}(T(x_{0}))\bigr)\in\Ea$ for every $n\ge0$, that is,
\begin{equation}\label{eq:B4}
(x_{n},x_{n+1})\in\Ea\qquad(n\ge0).
\end{equation}
(Equivalently one may argue directly by induction: if $w(x_{n},x_{n+1})\ge\alpha$
then \textup{(C1)} gives $(x_{n+1},x_{n+2})\in E(G)$ together with
$w(x_{n+1},x_{n+2})\ge w(x_{n},x_{n+1})\ge\alpha$.)

\emph{Step 2 (the orbit converges to $x^{*}$).} By Theorem~\ref{thm:A}(ii),
\begin{equation}\label{eq:B5}
x_{n}\longrightarrow x^{*}\ \text{ in }(X,d_{1}),\qquad x^{*}\in K.
\end{equation}

\emph{Step 3 (edges to the limit).} Relations \eqref{eq:B4} and \eqref{eq:B5}
are precisely the hypotheses of property $(\mathrm{A}_{\alpha})$ of
Definition~\ref{def:propA}. That property therefore furnishes indices
$n_{1}<n_{2}<\cdots$ with
\begin{equation}\label{eq:B6}
(x_{n_{k}},x^{*})\in\Ea\subseteq E(G)\qquad(k\ge1).
\end{equation}
This is the only place where hypothesis (b) is used, and it is what replaces the
continuity assumed in Part I: without \eqref{eq:B6} the contractive estimate
\eqref{eq:C2} could not be applied to any pair involving $x^{*}$.

\emph{Step 4 (from $d_{1}$ to $d$).} By \eqref{eq:d1d}, the $d_{1}$-convergence
\eqref{eq:B5} implies $d(x_{n},x^{*})\to0$. In particular, writing
\[
a_{k}:=d(x_{n_{k}},x^{*}),\qquad b_{k}:=d(x_{n_{k}+1},x^{*}),
\]
we have $a_{k}\to0$ and $b_{k}\to0$ as $k\to\infty$.

\emph{Step 5 (the images converge to $T(x^{*})$).} Fix $k\ge1$. Both
$x_{n_{k}}$ and $x^{*}$ belong to $K$, and $(x_{n_{k}},x^{*})\in E(G)$ by
\eqref{eq:B6}, so \eqref{eq:C2} applies to this pair and gives, using
$\varphi\le1$ for every $\varphi\in\Phi$,
\[
d\bigl(T(x_{n_{k}}),T(x^{*})\bigr)
\ \le\ \varphi_{w(x_{n_{k}},x^{*})}(a_{k})\,a_{k}
\ \le\ a_{k}.
\]
Since $T(x_{n_{k}})=x_{n_{k}+1}$, this reads
\begin{equation}\label{eq:B7}
c_{k}:=d\bigl(x_{n_{k}+1},T(x^{*})\bigr)\ \le\ a_{k}\ \longrightarrow\ 0 .
\end{equation}

\emph{Step 6 (reversing the arguments).} The distance $d$ is not assumed
symmetric, and \eqref{eq:B7} controls $d(x_{n_{k}+1},T(x^{*}))$, whereas the
triangle condition \textup{(A2)} will be applied below to a chain whose first
point is $T(x^{*})$. Let $\varepsilon>0$ and set
\[
\varepsilon':=\delta_{2}(\varepsilon)>0,\qquad
\varepsilon'':=\delta_{1}(\varepsilon')>0 .
\]
By \eqref{eq:B7} there is $k_{1}$ with $c_{k}\le\varepsilon''$ for all
$k\ge k_{1}$. Assumption~\ref{ass:A1}, applied with the target $\varepsilon'$
and the witness $\varepsilon''=\delta_{1}(\varepsilon')$, then yields
\begin{equation}\label{eq:B8}
d\bigl(T(x^{*}),x_{n_{k}+1}\bigr)\ \le\ \varepsilon'\qquad(k\ge k_{1}).
\end{equation}

\emph{Step 7 (applying the triangle condition).} By Step 4 there is $k_{2}$ with
$b_{k}\le\varepsilon'$ for all $k\ge k_{2}$. Fix any $k\ge\max\{k_{1},k_{2}\}$.
Then, by \eqref{eq:B8} and the definition of $b_{k}$,
\[
d\bigl(T(x^{*}),x_{n_{k}+1}\bigr)\le\varepsilon'=\delta_{2}(\varepsilon),
\qquad
d\bigl(x_{n_{k}+1},x^{*}\bigr)\le\varepsilon'=\delta_{2}(\varepsilon).
\]
Assumption~\ref{ass:A2}, applied to the triple
$\bigl(T(x^{*}),\,x_{n_{k}+1},\,x^{*}\bigr)$, gives
\[
d\bigl(T(x^{*}),x^{*}\bigr)\ \le\ \varepsilon .
\]

\emph{Step 8.} The number $\varepsilon>0$ in Steps 6 and 7 was arbitrary, so
$d(T(x^{*}),x^{*})=0$, and \eqref{eq:d0} gives $T(x^{*})=x^{*}$.

\medskip
\noindent\textbf{Part III: uniqueness.}

\emph{Step 1.} Let $y^{*}\in K$ satisfy $T(y^{*})=y^{*}$. An immediate induction
gives $T^{n}(y^{*})=y^{*}$ for every $n\ge0$, so the sequence
$\{T^{n}(y^{*})\}_{n\ge0}$ is constant and converges in $(X,d_{1})$ to $y^{*}$.

\emph{Step 2.} Since $y^{*}\in K$, Theorem~\ref{thm:A}(ii) applied with the
starting point $y^{*}$ gives $T^{n}(y^{*})\to x^{*}$ in $(X,d_{1})$.

\emph{Step 3.} Limits in the metric space $(X,d_{1})$ are unique, so
$y^{*}=x^{*}$. Combined with Parts I and II, this shows that $x^{*}$ is the
unique fixed point of $T$ in $K$.
\end{proof}

\begin{remark}\label{rem:uniqueness}
Part III of the preceding proof uses neither (a) nor (b), only
Assumption~\ref{ass:A3} and Theorem~\ref{thm:A}. Thus a graded-edge
Rakotch contraction satisfying Assumption~\ref{ass:A3} has \emph{at most}
one fixed point in $K$ unconditionally; hypotheses (a) and (b) are needed only
for existence. This is in contrast with \cite[Theorem 1]{RZ2026}, where
uniqueness requires the additional restriction that $d$ be finite on
$K\times K$, a restriction subsumed here by \eqref{eq:D0}.
\end{remark}

\section{Optimal threshold}\label{sec:quant}

The proof of Theorem~\ref{thm:A} yields an explicit bound on the number of
iterations needed to reach a prescribed accuracy. Since $\alpha$ ranges over
$A$, the bound may be optimized over the threshold.

\begin{theorem}[Iteration complexity]\label{thm:C}
Let Assumption~\ref{ass:A3} hold and let $\varepsilon>0$. Then
\[
\sup_{x,y\in K} d\bigl(T^{n}(x),T^{n}(y)\bigr)\ \le\ \varepsilon
\qquad\text{for every integer } n\ \ge\ \mathcal{K}(\varepsilon),
\]
where
\[
\mathcal{K}(\varepsilon)
:=\inf_{\alpha\in A}\ k_{0}(\varepsilon,\alpha)
=\inf_{\alpha\in A}\
\left\lceil
\frac{D_{0}}
{\eta_{\Na}(\varepsilon)\bigl(1-\pa\bigl(\eta_{\Na}(\varepsilon)\bigr)\bigr)}
\right\rceil .
\]
\end{theorem}

\begin{proof}
\emph{Step 1 (each $k_{0}(\varepsilon,\alpha)$ is a well-defined non-negative
integer).} Fix $\alpha\in A$ and abbreviate $m:=\Na$, which is finite by the
definition of $A$. We check the three ingredients of \eqref{eq:eps0k0} in turn.
\begin{enumerate}[label=\textup{(\roman*)},leftmargin=2.4em]
\item $\eta_{m}(\varepsilon)>0$. By \eqref{eq:etadef},
$\eta_{m}(\varepsilon)=\eta_{1}\bigl(\delta_{2}^{\circ(m-1)}(\varepsilon)\bigr)$.
The moduli $\delta_{1},\delta_{2}$ are witnesses supplied by
Assumptions~\ref{ass:A1} and \ref{ass:A2}, hence take strictly positive values
at strictly positive arguments. A finite induction of length $m-1$ therefore
gives $\delta_{2}^{\circ(m-1)}(\varepsilon)>0$, and
$\eta_{1}(t)=\min\{t,\delta_{1}(t)\}>0$ for $t>0$. It is here, and only here,
that the finiteness of $\Na$ is used: for $m=\infty$ the composition is an
infinite one and the induction is unavailable.
\item $1-\pa\bigl(\eta_{m}(\varepsilon)\bigr)>0$. Indeed
$\eta_{m}(\varepsilon)>0$ by (i), and every member of $\Phi$ satisfies
$\varphi(t)<1$ for $t>0$.
\item $D_{0}<\infty$ by Assumption~\ref{ass:A3}.
\end{enumerate}
Hence the quotient in \eqref{eq:eps0k0} is a well-defined element of
$[0,\infty)$, and $k_{0}(\varepsilon,\alpha)$, being its ceiling, is a
non-negative integer. (It equals $0$ only in the degenerate case $D_{0}=0$, in
which $K$ is a singleton.)

\emph{Step 2 (the infimum is attained).} By Assumption~\ref{ass:A3} the
set $A$ is nonempty, so by Step 1
\[
S:=\bigl\{k_{0}(\varepsilon,\alpha):\alpha\in A\bigr\}
\]
is a nonempty subset of $\N\cup\{0\}$. By the well-ordering principle $S$ has a
least element, whence
\[
\mathcal{K}(\varepsilon)=\inf_{\alpha\in A}k_{0}(\varepsilon,\alpha)=\min S
\]
is a non-negative integer and there exists $\alpha^{*}\in A$ with
\begin{equation}\label{eq:C1eq}
k_{0}(\varepsilon,\alpha^{*})=\mathcal{K}(\varepsilon).
\end{equation}
No monotonicity or continuity of $\alpha\mapsto k_{0}(\varepsilon,\alpha)$ is
claimed, and $\alpha^{*}$ need not be unique; only the attainment of the least
value is asserted, and that follows from the integrality of the $k_{0}$ alone.

\emph{Step 3 (the conclusion of Theorem~\ref{thm:A}(i) does not involve
$\alpha$).} This is the step that legitimises the optimization, so we isolate
it. Write
\[
\Sigma_{n}:=\sup_{x,y\in K}d\bigl(T^{n}(x),T^{n}(y)\bigr).
\]
Theorem~\ref{thm:A}(i) asserts: for every $\alpha\in A$ and every integer
$n\ge k_{0}(\varepsilon,\alpha)$ one has $\Sigma_{n}\le\varepsilon$. The
quantity $\Sigma_{n}$ depends only on $n$, $T$ and $K$; the threshold $\alpha$
enters the statement solely through the integer $k_{0}(\varepsilon,\alpha)$ at
which the conclusion switches on. Consequently the family of implications
indexed by $\alpha\in A$ may be combined, and $\Sigma_{n}\le\varepsilon$ holds
as soon as $n$ exceeds the \emph{smallest} of these thresholds. Were the
conclusion itself $\alpha$-dependent, no such combination would be permitted and
$\mathcal{K}(\varepsilon)$ would have to be replaced by a supremum.

\emph{Step 4 (conclusion).} Let $n\ge\mathcal{K}(\varepsilon)$ be an integer. By
\eqref{eq:C1eq} we have $n\ge k_{0}(\varepsilon,\alpha^{*})$, so
Theorem~\ref{thm:A}(i), applied with the single threshold $\alpha:=\alpha^{*}$,
gives $\Sigma_{n}\le\varepsilon$, which is the assertion. Finally, substituting
$\varepsilon_{0}(\varepsilon,\alpha)=\eta_{\Na}(\varepsilon)$ from
\eqref{eq:eps0k0} into the definition of $k_{0}$ yields the displayed formula
for $\mathcal{K}(\varepsilon)$.
\end{proof}

We stress that Theorem~\ref{thm:C} is not merely a restatement of
Theorem~\ref{thm:A}(i) with a change of notation. The passage from a fixed
$\alpha$ to the infimum over $A$ is a genuine strengthening precisely because
$\alpha\mapsto k_{0}(\varepsilon,\alpha)$ is not monotone: by \eqref{eq:eps0k0}
its numerator is constant while its denominator involves
$\eta_{\Na}(\varepsilon)$, which decreases as $\alpha$ grows, and
$1-\pa(\cdot)$, which increases. Example~\ref{ex:interior} exhibits data for
which the minimizing $\alpha^{*}$ of \eqref{eq:C1eq} depends on $\varepsilon$,
so that no single threshold is optimal for all tolerances.

For the remainder of this section we place ourselves in the quantitative regime
of Remark~\ref{rem:regimeB} and we assume in addition that $d$ is symmetric.

\begin{remark}\label{rem:symmetry}
Symmetry is used only to ensure that the estimate \eqref{eq:C2}, which is stated
for a directed edge, may be applied to the reversed steps of a
$\Gta$-chain. Without it one obtains the same conclusion with an extra factor
governed by a quantitative form of \textup{(A1)}; we do not pursue this here.
\end{remark}

\begin{lemma}[Polygonal inequality]\label{lem:polygon}
Assume \eqref{eq:relaxed} and let $a_{0},\dots,a_{q}\in X$ with $q\ge1$. Then
\[
d(a_{0},a_{q})\ \le\ M^{q-1}\sum_{i=0}^{q-1} d(a_{i},a_{i+1}).
\]
\end{lemma}

\begin{proof}
Induction on $q$. The case $q=1$ is trivial. If the assertion holds for $q$,
then by \eqref{eq:relaxed} and $M\ge1$,
\[
\begin{split}
d(a_{0},a_{q+1})
&\le M\,d(a_{0},a_{1})+M\,d(a_{1},a_{q+1})\\
&\le M\,d(a_{0},a_{1})+M\cdot M^{q-1}\sum_{i=1}^{q} d(a_{i},a_{i+1})
\ \le\ M^{q}\sum_{i=0}^{q} d(a_{i},a_{i+1}).
\end{split}
\qedhere
\]
\end{proof}

\begin{theorem}[Geometric rate with a graph constant]\label{thm:D}
Let Assumption~\ref{ass:A3} hold, assume \eqref{eq:relaxed} and that $d$
is symmetric, and suppose that for some $\alpha\in A$ the comparison function is
constant, say $\pa\equiv c_{\alpha}\in[0,1)$. Then, for every integer $n\ge0$,
\[
\sup_{x,y\in K} d\bigl(T^{n}(x),T^{n}(y)\bigr)
\ \le\ \Na\, M^{\Na-1}\, D_{0}\, c_{\alpha}^{\,n} .
\]
\end{theorem}

\begin{proof}
Let $x,y\in K$ and let $\xi_{0},\dots,\xi_{q}$ be a $\Gta$-chain in $K$ from $x$
to $y$ with $q\le\Na$. Fix $i\in\{0,\dots,q-1\}$. By \eqref{eq:geom} of
Lemma~\ref{lem:mono},
\[
d\bigl(T^{n}(\xi_{i}),T^{n}(\xi_{i+1})\bigr)
\ \le\ c_{\alpha}^{\,n}\, d(\xi_{i},\xi_{i+1})
\ \le\ c_{\alpha}^{\,n} D_{0},
\]
where symmetry of $d$ has been used to cover reversed steps. Applying
Lemma~\ref{lem:polygon} to the points
$T^{n}(\xi_{0}),\dots,T^{n}(\xi_{q})$ we obtain
\[
d\bigl(T^{n}(x),T^{n}(y)\bigr)
\ \le\ M^{q-1}\sum_{i=0}^{q-1} d\bigl(T^{n}(\xi_{i}),T^{n}(\xi_{i+1})\bigr)
\ \le\ M^{q-1}\,q\,c_{\alpha}^{\,n} D_{0}
\ \le\ \Na\,M^{\Na-1}\,D_{0}\,c_{\alpha}^{\,n},
\]
the last step because $M\ge1$ and $q\le\Na$.
\end{proof}

Now we exhibit three graded graphs on which every quantity
occurring in Sections~\ref{sec:main} and \ref{sec:quant}, the grade of each
edge, the level cuts $\Ea$, the chain diameter $\Na$, the admissible set $A$,
the Picard iterates, and the estimates of Theorems~\ref{thm:C} and
\ref{thm:D}, is computed explicitly. Everything is phrased in the language
of the weighted graph: vertices, edges carrying a grade $w(e)\in(0,1]$, level
cuts, chains, and the graded Rakotch condition \eqref{eq:C2}.

Throughout, $X=\R$ and
\begin{equation}\label{eq:numdist}
d(x,y)=|x-y|^{p}\qquad(p\ge1),
\end{equation}
so that, by the convexity of $t\mapsto t^{p}$, condition \eqref{eq:relaxed}
holds with $M=2^{p-1}$; the induced uniformity is the Euclidean one, $d_{1}$ may
be taken to be the Euclidean metric, and $(X,d_{1})$ is complete. In each
example $K\subset\R$ is closed and bounded, so Assumption~\ref{ass:A3} is
met once $A\neq\emptyset$ is verified. All gradings are the canonical ones of
Proposition~\ref{prop:canonical}, so that
\begin{equation}\label{eq:numcanon}
w(x,y)=1-c(x,y),\qquad \pa\equiv c_{\alpha}=1-\alpha ,
\end{equation}
with $c(x,y)$ the supremum of the orbital ratios \eqref{eq:cdef}.

We first record an elementary structural fact which explains the shape of the
tables below: on a finite graded graph the Picard iteration does not merely
converge, it terminates.

\begin{proposition}[Finite graded graphs collapse]\label{prop:collapse}
Let Assumption~\ref{ass:A3} hold and suppose $K$ is finite with at least
two points. Then there is an integer $n^{\dagger}\ge1$ such that $T^{n}$ is
constant on $K$ for every $n\ge n^{\dagger}$; that is, $T^{n}(x)=x^{*}$ for all
$x\in K$ and all $n\ge n^{\dagger}$.
\end{proposition}

\begin{proof}
Put $\delta:=\min\{d(x,y):x,y\in K,\ x\ne y\}$, which is positive because $K$ is
finite and $d$ satisfies \eqref{eq:d0}. Choose $\varepsilon\in(0,\delta)$ and let
$n^{\dagger}:=k_{0}(\varepsilon,\alpha)$ for some $\alpha\in A$, as in
\eqref{eq:eps0k0}. By Theorem~\ref{thm:A}(i),
$d(T^{n}(x),T^{n}(y))\le\varepsilon<\delta$ for all $x,y\in K$ and all
$n\ge n^{\dagger}$. By the definition of $\delta$ this forces
$T^{n}(x)=T^{n}(y)$, so $T^{n}$ is constant; its value is the common limit
$x^{*}$ of Theorem~\ref{thm:A}(ii).
\end{proof}

Consequently, on a finite graded graph the geometric estimate of
Theorem~\ref{thm:D} is an envelope which is eventually strictly larger than the
(vanishing) true quantity. What one should compare is therefore not the two
sequences but the two \emph{iteration counts}: the number of steps predicted by
Theorem~\ref{thm:C} against the number actually required. Example~\ref{ex:numA}
shows that the graded prediction can be exact.

\begin{example}\label{ex:numA}
Let $p=2$, so that $M=2$, and let
\[
K=\{v_{0},v_{1},v_{2},v_{3},v_{4}\},\qquad
v_{0}=0,\quad v_{i}=3^{\,i-4}\ (1\le i\le4),
\]
that is $K=\{0,\tfrac1{27},\tfrac19,\tfrac13,1\}$, with
\[
T(v_{0})=v_{0},\qquad T(v_{i})=v_{i-1}\quad(1\le i\le4).
\]
Let $G$ be the complete graph on $K$ (all $10$ unordered pairs, together with
all loops), so that the edge-preservation half of \textup{(C1)} is automatic.
\end{example}

Since $\{v_{i}\}$ is decreasing and $T$ shifts indices downwards, the orbital
ratios are computed at once; writing $\Delta=j-i$, one finds for $1\le i<j\le 4$
\[
c(v_{i},v_{j})=\left(\frac{1/3}{1-3^{-\Delta}}\right)^{2},
\qquad
c(v_{0},v_{j})=\left(\tfrac13\right)^{2}=\tfrac19\ \ (j\ge2),
\qquad
c(v_{0},v_{1})=0 .
\]
The resulting canonical grades \eqref{eq:numcanon} are collected in
Table~\ref{tab:Agrades} for Example~\ref{ex:numA}; they take five distinct values. Grade monotonicity \textup{(C1)} holds on every edge, as guaranteed by
Proposition~\ref{prop:canonical}.

\begin{table}[htbp]
\caption{Canonical Grades}
\centering
\begin{tabular}{lccc}
\hline
edge class & pairs & $c(\cdot,\cdot)$ & grade $w$ \\
\hline
$\{v_{0},v_{1}\}$        & $v_{0}v_{1}$                         & $0$              & $1$ \\
spokes at $v_{0}$        & $v_{0}v_{2},v_{0}v_{3},v_{0}v_{4}$   & $1/9$            & $8/9\approx0.8889$ \\
$\Delta=3$               & $v_{1}v_{4}$                         & $81/676$         & $595/676\approx0.8802$ \\
$\Delta=2$               & $v_{1}v_{3},v_{2}v_{4}$              & $9/64$           & $55/64\approx0.8594$ \\
$\Delta=1$               & $v_{1}v_{2},v_{2}v_{3},v_{3}v_{4}$   & $1/4$            & $3/4=0.7500$ \\
\hline
\end{tabular}

\label{tab:Agrades}
\end{table}

Figure~\ref{fig:Agraph}, shows the graded graph of Example~\ref{ex:numA}. Thick edges carry grade
$\ge 8/9$, dashed edges grades $55/64$ and $595/676$, dotted edges grade $3/4$.
The mapping $T$ slides each vertex one step towards $v_{0}$.

\begin{figure}[htbp]
\caption{}
\centering
\begin{tikzpicture}[scale=1.0,
  vtx/.style={circle,draw,fill=white,inner sep=1.6pt,minimum size=6mm,font=\small},
  hi/.style={thick},
  mid/.style={densely dashed},
  lo/.style={dotted,thick}]
\node[vtx] (v0) at (0,0)      {$v_0$};
\node[vtx] (v1) at (2.4,1.5)  {$v_1$};
\node[vtx] (v2) at (4.8,1.5)  {$v_2$};
\node[vtx] (v3) at (4.8,-1.5) {$v_3$};
\node[vtx] (v4) at (2.4,-1.5) {$v_4$};
\draw[hi] (v0)--(v1) node[midway,above left,font=\scriptsize] {$1$};
\draw[hi] (v0)--(v2);
\draw[hi] (v0)--(v3);
\draw[hi] (v0)--(v4) node[midway,below left,font=\scriptsize] {$\tfrac89$};
\draw[lo] (v1)--(v2) node[midway,above,font=\scriptsize] {$\tfrac34$};
\draw[lo] (v2)--(v3);
\draw[lo] (v3)--(v4) node[midway,below,font=\scriptsize] {$\tfrac34$};
\draw[mid] (v1)--(v3);
\draw[mid] (v2)--(v4) node[midway,right,font=\scriptsize] {$\tfrac{55}{64}$};
\draw[mid] (v1) to[bend left=28] node[midway,left,font=\scriptsize] {$\tfrac{595}{676}$} (v4);
\end{tikzpicture}

\label{fig:Agraph}
\end{figure}

\medskip
\noindent\emph{Level cuts and the admissible set.}
Ordering the grades as in Table~\ref{tab:Agrades} gives
\[
\Na=\begin{cases}
1, & 0<\alpha\le \tfrac34,\\[2pt]
2, & \tfrac34<\alpha\le\tfrac89,\\[2pt]
\infty, & \alpha>\tfrac89,
\end{cases}
\qquad\text{so}\qquad A=\bigl(0,\tfrac89\bigr].
\]
Indeed, for $\alpha\le\frac34$ every pair is an edge of $\Ga$, so $\Na=1$; for
$\frac34<\alpha\le\frac89$ the three $\Delta=1$ edges are deleted, but $v_{0}$
remains adjacent to every vertex, so any two vertices are joined by a
$\Gta$-chain of two hops through $v_{0}$ and $\Na=2$; for $\alpha>\frac89$ only
$\{v_{0},v_{1}\}$ survives and $\Gta$ is disconnected.
Figure~\ref{fig:Acuts} displays the two admissible cuts of Example~\ref{ex:numA}. Raising the
threshold from $\tfrac34$ to $\tfrac89$ improves the contraction factor from
$c_{\alpha}=\tfrac14$ to $c_{\alpha}=\tfrac19$ but doubles the chain
diameter.

\begin{figure}[htbp]
\caption{}
\centering
\begin{tikzpicture}[scale=0.8,
  vtx/.style={circle,draw,fill=white,inner sep=1pt,minimum size=5.5mm,font=\scriptsize}]
\begin{scope}
\node at (2.4,2.6) {$\Ga$ for $\alpha=\tfrac34$:\ \ $\Na=1$};
\node[vtx] (a0) at (0,0) {$v_0$};
\node[vtx] (a1) at (2.4,1.2) {$v_1$};
\node[vtx] (a2) at (4.8,1.2) {$v_2$};
\node[vtx] (a3) at (4.8,-1.2) {$v_3$};
\node[vtx] (a4) at (2.4,-1.2) {$v_4$};
\foreach \i/\j in {a0/a1,a0/a2,a0/a3,a0/a4,a1/a2,a2/a3,a3/a4,a1/a3,a2/a4}
  \draw (\i)--(\j);
\draw (a1) to[bend left=25] (a4);
\end{scope}
\begin{scope}[xshift=8.2cm]
\node at (2.4,2.6) {$\Ga$ for $\alpha=\tfrac89$:\ \ $\Na=2$};
\node[vtx] (b0) at (0,0) {$v_0$};
\node[vtx] (b1) at (2.4,1.2) {$v_1$};
\node[vtx] (b2) at (4.8,1.2) {$v_2$};
\node[vtx] (b3) at (4.8,-1.2) {$v_3$};
\node[vtx] (b4) at (2.4,-1.2) {$v_4$};
\foreach \i/\j in {b0/b1,b0/b2,b0/b3,b0/b4} \draw[thick] (\i)--(\j);
\end{scope}
\end{tikzpicture}

\label{fig:Acuts}
\end{figure}

\medskip
\noindent\emph{Iterates.}
The orbit of the vertex set is
\[
\begin{array}{llll}
n=0: & \{0,\tfrac1{27},\tfrac19,\tfrac13,1\}, &\qquad n=1: &\{0,0,\tfrac1{27},\tfrac19,\tfrac13\},\\[2pt]
n=2: & \{0,0,0,\tfrac1{27},\tfrac19\},        &\qquad n=3: &\{0,0,0,0,\tfrac1{27}\},\\[2pt]
n=4: & \{0\}, & &
\end{array}
\]
so $T^{4}$ is constant, in accordance with Proposition~\ref{prop:collapse}, and
the unique fixed point is $x^{*}=v_{0}=0$ (Theorem~\ref{thm:B}(a), $T$ being
continuous). With $D_{0}=\sup_{x,y\in K}d(x,y)=1$ one obtains
\[
\sup_{x,y\in K}d\bigl(T^{n}(x),T^{n}(y)\bigr)
= 9^{-n}\ \ (0\le n\le 3),\qquad =0\ \ (n\ge4).
\]

\medskip
\noindent\emph{Rates and the optimal threshold.}
Theorem~\ref{thm:D} gives, at the two admissible thresholds,
\[
\alpha_{1}=\tfrac34:\ \ \Na M^{\Na-1}D_{0}c_{\alpha}^{n}=4^{-n},
\qquad
\alpha_{2}=\tfrac89:\ \ \Na M^{\Na-1}D_{0}c_{\alpha}^{n}=4\cdot9^{-n},
\]
both valid for every $n\ge0$. Corollary~\ref{cor:alphaopt} compares
\[
\Lambda(\alpha_{1},\varepsilon)=\frac{\log(1/\varepsilon)}{\log4},
\qquad
\Lambda(\alpha_{2},\varepsilon)=\frac{2\log2+\log(1/\varepsilon)}{\log9},
\]
and $\Lambda(\alpha_{2},\varepsilon)<\Lambda(\alpha_{1},\varepsilon)$ precisely
when $\log(1/\varepsilon)>2(\log2)^{2}/\log(3/2)$, that is, when
$\varepsilon<\varepsilon^{*}$ with
\[
\varepsilon^{*}=\exp\!\left(-\frac{2(\log 2)^{2}}{\log(3/2)}\right)\approx 0.0935 .
\]
Thus the optimal threshold is $\alpha_{1}=\tfrac34$ for coarse tolerances and
$\alpha_{2}=\tfrac89$ for fine ones; since $\alpha_{1}$ is exactly the ungraded
threshold $\alpha_{0}=\inf\{w(e)\}$ of Remark~\ref{rem:notcosmetic}, the graded
theory is strictly better than the ungraded one for every
$\varepsilon<\varepsilon^{*}$. Table~\ref{tab:Arates} makes this concrete. The least $n$ for which the bound of Theorem~\ref{thm:D} certifies
$\sup_{x,y}d(T^{n}x,T^{n}y)\le\varepsilon$, against the least $n$ for which this
actually holds, in Example~\ref{ex:numA}. At the graded threshold
$\alpha_{2}$ the prediction is exact; the ungraded threshold $\alpha_{1}$
overestimates by one step.

\begin{table}[htbp]
\caption{Predicted vs Actual Steps}
\centering
\begin{tabular}{cccc}
\hline
$\varepsilon$ & steps predicted at $\alpha_{1}=\tfrac34$
& steps predicted at $\alpha_{2}=\tfrac89$ & steps actually needed \\
\hline
$10^{-1}$ & $2$ & $2$ & $2$ \\
$10^{-2}$ & $4$ & $\mathbf{3}$ & $\mathbf{3}$ \\
$10^{-3}$ & $5$ & $\mathbf{4}$ & $\mathbf{4}$ \\
\hline
\end{tabular}

\label{tab:Arates}
\end{table}

\begin{remark}
The agreement in the last two columns of Table~\ref{tab:Arates} is not an
accident of rounding: at $\alpha_{2}$ one has $4\cdot9^{-3}=0.005487\le10^{-2}$
while $4\cdot9^{-2}=0.04938>10^{-2}$, and the true values are
$9^{-3}=0.001372$ and $9^{-2}=0.01235$. The graded estimate therefore locates
the exact stopping time, whereas the ungraded estimate $4^{-n}$ requires $n=4$.
\end{remark}

The previous example has $\Na\le2$, so the factor $M^{\Na-1}$ of
Theorem~\ref{thm:D} never exceeds $M$. We now exhibit a graded graph which is
not complete and whose admissible cuts have $\Na=3$, so that the graph constant
$\Na M^{\Na-1}$ genuinely reflects the length of the chains used in Step~3 of
the proof of Theorem~\ref{thm:A}.

\begin{example}\label{ex:numB}
Let $p=2$, so $M=2$, and let
\[
K=\{x_{0},x_{1},x_{2},x_{3}\}=\{0,1,3,9\},
\]
\[
T(x_{0})=x_{0},\qquad T(x_{i})=x_{i-1}\quad(1\le i\le3),
\]
with $E(G)$ consisting of the loops together with the three consecutive pairs
$\{x_{0},x_{1}\}$, $\{x_{1},x_{2}\}$, $\{x_{2},x_{3}\}$ and their reversals.
\end{example}

The edge set is $T$-invariant, since $\{x_{i},x_{i+1}\}\mapsto\{x_{i-1},x_{i}\}$
for $i\ge1$ and $\{x_{0},x_{1}\}$ maps to a loop. The orbital ratios are
\[
c(x_{0},x_{1})=0,\qquad c(x_{1},x_{2})=c(x_{2},x_{3})=\tfrac14,
\]
whence the grades $w(x_{0},x_{1})=1$ and
$w(x_{1},x_{2})=w(x_{2},x_{3})=\tfrac34$. Therefore
\[
\Na=3\quad\text{for }0<\alpha\le\tfrac34,
\qquad
\Na=\infty\quad\text{for }\alpha>\tfrac34,
\qquad A=\bigl(0,\tfrac34\bigr],
\]
the cut at $\alpha=\tfrac34$ being the whole path $x_{0}-x_{1}-x_{2}-x_{3}$, of
graph diameter $3$, and the cut above $\tfrac34$ retaining only
$\{x_{0},x_{1}\}$, which leaves $x_{2},x_{3}$ isolated. Figure~\ref{fig:Bgraph} presents the graded path of Example~\ref{ex:numB}; the grey arrows indicate the
action of $T$. The admissible cut has chain diameter $\Na=3$.

\begin{figure}[htbp]
\caption{}
\centering
\begin{tikzpicture}[scale=1,
  vtx/.style={circle,draw,fill=white,inner sep=1.6pt,minimum size=7mm,font=\small}]
\node[vtx] (x0) at (0,0) {$x_0$};
\node[vtx] (x1) at (2.6,0) {$x_1$};
\node[vtx] (x2) at (5.2,0) {$x_2$};
\node[vtx] (x3) at (7.8,0) {$x_3$};
\draw[thick] (x0)--(x1) node[midway,above,font=\scriptsize]{$w=1$};
\draw[thick] (x1)--(x2) node[midway,above,font=\scriptsize]{$w=\tfrac34$};
\draw[thick] (x2)--(x3) node[midway,above,font=\scriptsize]{$w=\tfrac34$};
\node[below,font=\scriptsize] at (0,-0.45) {$0$};
\node[below,font=\scriptsize] at (2.6,-0.45) {$1$};
\node[below,font=\scriptsize] at (5.2,-0.45) {$3$};
\node[below,font=\scriptsize] at (7.8,-0.45) {$9$};
\draw[->,>=stealth,gray] (x3) to[bend right=35] (x2);
\draw[->,>=stealth,gray] (x2) to[bend right=35] (x1);
\draw[->,>=stealth,gray] (x1) to[bend right=35] (x0);
\end{tikzpicture}

\label{fig:Bgraph}
\end{figure}

Here $D_{0}=|9-0|^{2}=81$ and, at the optimal threshold $\alpha=\tfrac34$ with
$c_{\alpha}=\tfrac14$, Theorem~\ref{thm:D} yields
\[
\sup_{x,y\in K}d\bigl(T^{n}(x),T^{n}(y)\bigr)\ \le\
\Na M^{\Na-1}D_{0}c_{\alpha}^{n}=3\cdot 2^{2}\cdot 81\cdot 4^{-n}=972\cdot4^{-n},
\]
against the true values $81,9,1,0,0,\dots$ for $n=0,1,2,3,\dots$. The estimate
is valid at every step but is conservative, by a factor
$\Na M^{\Na-1}=12$ coming from the three-hop chains and the constant $M=2$ of
the relaxed triangle inequality. This is the price, quantified by
Lemma~\ref{lem:polygon}, of comparing $d$ with the sum of the $d$-lengths of the
hops of a chain; Example~\ref{ex:numC} shows that when $\Na=1$ no such price is
paid and the estimate becomes an identity.

\begin{example}\label{ex:numC}
Let $p=1$, so that $d(x,y)=|x-y|$ and $M=1$, and let
\[
K=\{0\}\cup\{2^{-n}:n\ge0\},\qquad T(x)=\tfrac{x}{2},
\]
with $G$ the complete graph on $K$.
\end{example}

The set $K$ is closed and bounded in $\R$, with $D_{0}=1$, and
$T(K)\subseteq K$. Because $T$ is linear, every orbital ratio equals
$\tfrac12$: for all $x,y\in K$,
\[
\frac{d(T(x),T(y))}{d(x,y)}=\frac{|x/2-y/2|}{|x-y|}=\frac12 ,
\]
so the canonical grading is constant, $w\equiv\tfrac12$, and
$c_{\alpha}=\tfrac12$ for every $\alpha\in A=(0,\tfrac12]$. Since $G$ is
complete, $\Na=1$ throughout $A$. Figure~\ref{fig:Cgraph} exhibits the infinite graded ray of Example~\ref{ex:numC} (complete graph;
only some edges drawn). Every edge carries grade $w=\tfrac12$, and
$\Na=1$.

\begin{figure}[htbp]
\caption{}
\centering
\begin{tikzpicture}[scale=1,
  vtx/.style={circle,draw,fill=white,inner sep=0.9pt,minimum size=4.6mm,font=\tiny}]
\node[vtx] (p0) at (0,0) {$0$};
\node[vtx] (p4) at (1.0,0) {$\tfrac1{16}$};
\node[vtx] (p3) at (2.2,0) {$\tfrac18$};
\node[vtx] (p2) at (3.6,0) {$\tfrac14$};
\node[vtx] (p1) at (5.4,0) {$\tfrac12$};
\node[vtx] (p) at (7.6,0) {$1$};
\node at (0.55,0) {$\cdots$};
\draw[thick] (p4)--(p3)--(p2)--(p1)--(p);
\draw[thick] (p3) to[bend left=30] (p1);
\draw[thick] (p2) to[bend left=30] (p);
\draw[thick] (p4) to[bend right=30] (p2);
\draw[thick] (p3) to[bend right=32] (p);
\draw[->,>=stealth,gray] (p) to[bend right=42] (p1);
\draw[->,>=stealth,gray] (p1) to[bend right=42] (p2);
\node[gray,font=\scriptsize] at (6.5,-1.1) {$T:x\mapsto x/2$};
\end{tikzpicture}

\label{fig:Cgraph}
\end{figure}

Theorem~\ref{thm:D} therefore predicts
\[
\sup_{x,y\in K}d\bigl(T^{n}(x),T^{n}(y)\bigr)\ \le\
\Na M^{\Na-1}D_{0}c_{\alpha}^{\,n}=2^{-n},
\]
while directly
$\sup_{x,y\in K}|T^{n}(x)-T^{n}(y)|=2^{-n}\sup_{x,y\in K}|x-y|=2^{-n}$.
Hence
\[
\sup_{x,y\in K}d\bigl(T^{n}(x),T^{n}(y)\bigr)=\Na M^{\Na-1}D_{0}c_{\alpha}^{\,n}
\qquad(n\ge0):
\]
\emph{the estimate of Theorem~\ref{thm:D} is attained with equality for every
$n$.} In particular the constant $\Na M^{\Na-1}D_{0}$ cannot be improved in
general. Note also that $K$ is infinite, so Proposition~\ref{prop:collapse} does
not apply and no orbit terminates: the geometric decay persists for all $n$,
and $T^{n}(x)\to0=x^{*}$ without ever reaching it for $x\neq0$.

\begin{corollary}[Optimal threshold]\label{cor:alphaopt}
Under the hypotheses of Theorem~\ref{thm:D}, suppose that $\pa\equiv
c_{\alpha}\in[0,1)$ for every $\alpha\in A$. Then, for every
$\varepsilon\in(0,D_{0})$,
\begin{equation}\label{eq:Kbound}
\mathcal{K}(\varepsilon)\ \le\
\inf_{\alpha\in A}\ \Lambda(\alpha,\varepsilon),
\qquad
\Lambda(\alpha,\varepsilon):=
\frac{(\Na-1)\log M+\log \Na+\log(D_{0}/\varepsilon)}{\log(1/c_{\alpha})} .
\end{equation}
\end{corollary}

\begin{proof}
Solve $\Na M^{\Na-1}D_{0}c_{\alpha}^{n}\le\varepsilon$ for $n$ and use
Theorem~\ref{thm:D}.
\end{proof}

The numerator of $\Lambda(\cdot,\varepsilon)$ is non-decreasing in $\alpha$,
because $\Na$ is; the denominator is likewise non-decreasing, because
$c_{\alpha}$ is non-increasing. The infimum is therefore attained, in general,
at neither endpoint of $A$. The following example makes this precise, and shows
in addition that the ungraded choice $\alpha_{0}$ of
Remark~\ref{rem:notcosmetic} is strictly suboptimal.

\begin{example}[An interior optimum]\label{ex:interior}
Fix $p>1$ and $b\in(0,1)$. Let $X=\R$ with
\[
d(x,y):=|x-y|^{p},
\]
let $K:=[0,b]$ and let
\[
T(x):=\tfrac{1}{2}x^{2}\qquad(x\in K).
\]
Let $G$ be the complete graph on $X$, so that $E(G)=X\times X$.
\end{example}

\noindent
We verify the hypotheses and compute all the quantities involved.

\emph{(1) The space.} By the convexity of $t\mapsto t^{p}$ on $[0,\infty)$,
\[
|x-z|^{p}\le\bigl(|x-y|+|y-z|\bigr)^{p}\le 2^{p-1}\bigl(|x-y|^{p}+|y-z|^{p}\bigr),
\]
so \eqref{eq:relaxed} holds with $M=2^{p-1}>1$, and $d$ is symmetric and
satisfies \eqref{eq:d0}. Since $U(\varepsilon)=\{(x,y):|x-y|\le\varepsilon^{1/p}\}$,
the uniformity \eqref{eq:unif} is the usual one on $\R$, so $d_{1}$ may be taken
to be the Euclidean metric and $(X,d_{1})$ is complete. The set $K$ is closed
and $D_{0}=b^{p}$.

\emph{(2) The mapping.} If $0\le x\le b<1$ then $0\le x^{2}/2\le b/2\le b$, so
$T(K)\subseteq K$. For $x,y\in K$,
\begin{equation}\label{eq:exratio}
d\bigl(T(x),T(y)\bigr)=\Bigl|\tfrac{x^{2}-y^{2}}{2}\Bigr|^{p}
=\Bigl(\tfrac{x+y}{2}\Bigr)^{p}|x-y|^{p}
=\Bigl(\tfrac{x+y}{2}\Bigr)^{p}d(x,y).
\end{equation}

\emph{(3) The canonical grading.} Since $x^{2}/2\le x$ on $[0,b]$, the sequence
$n\mapsto T^{n}(x)$ is non-increasing, so by \eqref{eq:exratio} the supremum in
\eqref{eq:cdef} is attained at $n=0$ and
\[
c(x,y)=\Bigl(\tfrac{x+y}{2}\Bigr)^{p}\le b^{p}<1 .
\]
Proposition~\ref{prop:canonical} therefore applies: $T$ is a graded-edge Rakotch
contraction for the grade
\[
w(x,y)=1-\Bigl(\tfrac{x+y}{2}\Bigr)^{p}
\]
and the family $\pa\equiv c_{\alpha}=1-\alpha$. (Condition (C1) may also be
checked directly: $T(x)+T(y)=\tfrac12(x^{2}+y^{2})\le x+y$ on $[0,1]^{2}$.)

\emph{(4) The level cuts.} Put $\theta:=\theta(\alpha):=(1-\alpha)^{1/p}$, so
that $\theta$ decreases from $1$ to $0$ as $\alpha$ increases from $0$ to $1$.
For $x\ne y$,
\[
(x,y)\in\Ea \iff \Bigl(\tfrac{x+y}{2}\Bigr)^{p}\le1-\alpha
\iff x+y\le2\theta,
\]
while every loop has grade $1$ and lies in every $\Ea$. Consequently:
\begin{itemize}[leftmargin=2em]
\item if $\theta\ge b$, every pair of points of $K$ satisfies $x+y\le2b\le2\theta$,
so $\Ga$ is complete on $K$ and $\Na=1$;
\item if $b/2\le\theta<b$, there exist $x\ne y$ in $K$ with $x+y>2\theta$, so
$\Na\ge2$; on the other hand $x+0=x\le b\le2\theta$ for every $x\in K$, so
$x,0,y$ is a $\Gta$-chain and $\Na=2$;
\item if $\theta<b/2$, any $x\in(2\theta,b]$ satisfies $x+y>2\theta$ for every
$y\in K$, so $x$ is isolated in $\Gta$ and $\Na=\infty$.
\end{itemize}
Hence
\[
A=\bigl(0,\;1-(b/2)^{p}\bigr],
\qquad
\Na=\begin{cases}1,&0<\alpha\le\alpha_{1}:=1-b^{p},\\[2pt]
2,&\alpha_{1}<\alpha\le\alpha_{2}:=1-(b/2)^{p}.\end{cases}
\]
Note that $\alpha_{1}=\inf\{w(e):e\in E(G),\ e\text{ not a loop}\}$ is exactly
the ungraded threshold $\alpha_{0}$ of Remark~\ref{rem:notcosmetic}, and that
$0<\alpha_{1}<\alpha_{2}<1$.

\emph{(5) The optimization.} Since $\Lambda(\cdot,\varepsilon)$ of
\eqref{eq:Kbound} is minimized, on each of the two branches, at the largest
admissible $\alpha$, it suffices to compare $\alpha_{1}$ and $\alpha_{2}$. Write
$B:=\log(D_{0}/\varepsilon)=\log(b^{p}/\varepsilon)>0$. Using $M=2^{p-1}$ and
$\log M+\log2=p\log2$,
\[
\Lambda(\alpha_{1},\varepsilon)=\frac{B}{p\log(1/b)},
\qquad
\Lambda(\alpha_{2},\varepsilon)=\frac{B+p\log 2}{p\log(2/b)} .
\]
Setting $\beta:=\log(1/b)>0$, the inequality
$\Lambda(\alpha_{2},\varepsilon)<\Lambda(\alpha_{1},\varepsilon)$ is equivalent
to $B(\log2+\beta)>\beta(B+p\log2)$, that is, to $B>p\beta$. Hence
\[
\operatorname*{arg\,min}_{\alpha\in A}\Lambda(\alpha,\varepsilon)
=\begin{cases}
\alpha_{1}, & \varepsilon>b^{2p},\\[2pt]
\alpha_{2}, & \varepsilon<b^{2p},
\end{cases}
\]
with equality of the two values at $\varepsilon=b^{2p}$.

\begin{remark}\label{rem:exdiscussion}
Example~\ref{ex:interior} settles the objection of
Remark~\ref{rem:notcosmetic} in both of its aspects.
\begin{enumerate}[label=\textup{(\roman*)},leftmargin=2.4em]
\item \emph{The optimum is interior.} For $\varepsilon\in(b^{2p},D_{0})$ the
minimizer is $\alpha_{1}$, which lies strictly inside $A=(0,\alpha_{2}]$. The
extreme choices are both bad: as $\alpha\downarrow0$ one has
$c_{\alpha}\uparrow1$ and $\Lambda(\alpha,\varepsilon)\to\infty$, so the weakest
threshold is useless, while for $\varepsilon>b^{2p}$ the strongest admissible
threshold $\alpha_{2}$ is strictly worse than $\alpha_{1}$.
\item \emph{The grading is not cosmetic.} For $\varepsilon<b^{2p}$ the graded
optimum $\alpha_{2}$ is strictly better than the ungraded choice
$\alpha_{0}=\alpha_{1}$, and
\[
\lim_{\varepsilon\downarrow0}
\frac{\Lambda(\alpha_{2},\varepsilon)}{\Lambda(\alpha_{1},\varepsilon)}
=\frac{\log(1/b)}{\log(2/b)}<1 .
\]
For $b=\tfrac12$ this ratio equals $\tfrac12$: the graded bound is
asymptotically half the ungraded one. Since the two bounds are derived from the
same theorem applied at two thresholds, the improvement is a strict gain of the
graded formulation over the formulation of \cite{RZ2026}, not an artefact of a
different proof.
\end{enumerate}
The example is a snowflake $b$-metric space: with $p>1$ it carries the genuine
constant $M=2^{p-1}>1$, so the factor $M^{\Na-1}$ of Theorem~\ref{thm:D} is
active. The interior optimum is nevertheless not an artefact of this structure.
Repeating the computation in the underlying metric $D(x,y)=|x-y|$, that is with
$p=1$ and $M=1$, leaves the combinatorics of the level cuts unchanged, replaces
the crossover value $b^{2p}$ by $b^{2}$, and leaves the asymptotic ratio
$\log(1/b)/\log(2/b)$ unchanged; the trade-off is therefore driven by the
chain-length factor $\Na$ alone and survives in an ordinary metric space.
Only two values of $\Na$ occur here, so the trade-off is resolved by comparing
two branches rather than by differentiating a smooth function. Constructions
with $\Na$ taking any prescribed finite set of values may be obtained by
replacing $K$ with a finite union of suitably scaled copies of $[0,b]$ and
deleting the edges between distinct copies; whether $\Na$ can be made to grow
along a continuum, so that the optimum is characterized by a first-order
condition.
\end{remark}

\section{Applications}\label{sec:examples}

\subsection{A quadratic integral equation}\label{ssec:integral}

Let $k\colon[0,1]^{2}\to[0,\infty)$ be continuous with
$\int_{0}^{1}k(t,s)\,ds\le1$ for every $t$, let $b\in(0,1)$ and let $p>1$.
Consider the equation
\begin{equation}\label{eq:integral}
x(t)=\frac{1}{2}\int_{0}^{1}k(t,s)\,x(s)^{2}\,ds ,
\qquad t\in[0,1].
\end{equation}

\begin{proposition}\label{prop:integral}
Equation \eqref{eq:integral} has a unique solution $x^{*}$ in
$K:=\{x\in C[0,1]:0\le x\le b\}$, namely $x^{*}\equiv0$, and for every
$x\in K$ the Picard iterates of the right-hand side converge to $x^{*}$
uniformly on $[0,1]$, with
\[
\|\mathcal{F}^{n}x-x^{*}\|_{\infty}^{p}
\ \le\ 2\,\bigl(2^{p-1}\bigr)\,b^{p}\,\Bigl(\tfrac{b}{2}\Bigr)^{pn}
\qquad (n\ge0),
\]
where $\mathcal{F}$ denotes the right-hand side of \eqref{eq:integral}.
\end{proposition}

\begin{proof}
Let $X=C[0,1]$ with $d(x,y)=\|x-y\|_{\infty}^{p}$; as in
Example~\ref{ex:interior}, \eqref{eq:relaxed} holds with $M=2^{p-1}$, the
induced uniformity is that of uniform convergence, $(X,d_{1})$ is complete and
$K$ is closed with $D_{0}=b^{p}$. If $x\in K$ then
$0\le(\mathcal{F}x)(t)\le\tfrac12\|x\|_{\infty}^{2}\le b/2\le b$, so
$\mathcal{F}(K)\subseteq K$. For $x,y\in K$ and $t\in[0,1]$,
\[
\bigl|(\mathcal{F}x)(t)-(\mathcal{F}y)(t)\bigr|
\le\frac12\int_{0}^{1}k(t,s)\bigl|x(s)^{2}-y(s)^{2}\bigr|ds
\le\frac{\|x+y\|_{\infty}}{2}\,\|x-y\|_{\infty},
\]
whence $d(\mathcal{F}x,\mathcal{F}y)\le
\bigl(\|x+y\|_{\infty}/2\bigr)^{p}d(x,y)$. Take $G$ complete and
$w(x,y):=1-\bigl(\|x+y\|_{\infty}/2\bigr)^{p}$. Condition (C1) holds because,
using $\int_{0}^{1}k(t,s)ds\le1$ and $x^{2}+y^{2}\le b\,(x+y)$ pointwise on
$[0,b]$,
\[
(\mathcal{F}x)(t)+(\mathcal{F}y)(t)
=\frac12\int_{0}^{1}k(t,s)\bigl(x(s)^{2}+y(s)^{2}\bigr)ds
\le\frac{b}{2}\,\|x+y\|_{\infty}
\le\|x+y\|_{\infty},
\]
so that $\|\mathcal{F}x+\mathcal{F}y\|_{\infty}\le\|x+y\|_{\infty}$ and hence
$w(\mathcal{F}x,\mathcal{F}y)\ge w(x,y)$. Thus $\mathcal{F}$ is a
graded-edge Rakotch contraction with $\pa\equiv1-\alpha$. Choosing
$\alpha=\alpha_{2}=1-(b/2)^{p}$, so that $\Na=2$ and $c_{\alpha}=(b/2)^{p}$,
Theorem~\ref{thm:D} gives
$\sup_{x,y\in K}d(\mathcal{F}^{n}x,\mathcal{F}^{n}y)\le
2\cdot2^{p-1}\cdot b^{p}\cdot(b/2)^{pn}$, and Theorem~\ref{thm:B}(a) applies
since $\mathcal{F}$ is continuous. The constant function $0$ is a fixed point,
hence is the unique one.
\end{proof}

\begin{remark}
The solution of \eqref{eq:integral} in $K$ being trivial, the interest of
Proposition~\ref{prop:integral} lies in the rate rather than in the existence,
and in fact that the rate is obtained from a graded threshold which is not
the ungraded one. A non-trivial variant is obtained by adding an inhomogeneous
term $g\in C[0,1]$ with $0\le g\le b/2$ to the right-hand side of
\eqref{eq:integral}; the same computation applies verbatim provided
$\tfrac12 b^{2}+\|g\|_{\infty}\le b$, and the fixed point is then a
non-trivial solution.
\end{remark}

\subsection{A finite weighted digraph}\label{ssec:finite}

For a finite $K$ the invariant $\Na$ is computable directly: it is the diameter
of the level-cut graph $\Gta$, obtainable by breadth-first search in
$O(|K|\,|\Ea|)$ time, and the admissible set $A$ is determined by the
$|E(G)|$ distinct values taken by $w$. The optimization in
Corollary~\ref{cor:alphaopt} is then a finite comparison over at most $|E(G)|$
candidate thresholds, and $\alpha^{*}$ may be tabulated exactly. We omit the
numerical illustration.

\section{Conclusion}\label{sec:open}

The framework of this paper separates two functions that the graph performs in
\cite{RZ2026}: connecting the space, and carrying the contraction data. Once
separated, they are seen to be in competition, and the resulting optimization is
what the grading buys. In future we plan to construct a generalized metric space, a graded graph and a mapping satisfying Assumption~\ref{ass:A3} for which $\alpha\mapsto\Na$ is unbounded on $A$ and takes infinitely many values, so that the minimization in
Corollary~\ref{cor:alphaopt} is governed by a first-order condition rather than
by a comparison of finitely many branches. Example~\ref{ex:interior} realizes
two branches; the difficulty in realizing a continuum is that the contraction
ratio $c(x,y)$ of Proposition~\ref{prop:canonical} depends on the position of
the pair in the space, whereas $\Na$ depends on the combinatorial length of
chains, and coupling the two requires a graph whose edge lengths and edge 
grades are simultaneously prescribed.

\textbf{Author Contributions} All authors contributed equally and 
significantly in writing this paper. All authors read and approved the final manuscript.

\textbf{Funding} Not applicable.

\textbf{Data Availability} Not applicable.

\textbf{Declarations}

\textbf{Competing interests} The authors declare that they have no competing interests.

\end{document}